\documentclass{amsart}

\newtheorem{theorem}{Theorem}[section]

\newtheorem{lemma}[theorem]{Lemma}
\newtheorem{proposition}[theorem]{Proposition}
\newtheorem{definition}[theorem]{Definition}
\newtheorem{example}[theorem]{Example}
\newtheorem{conjecture}[theorem]{Conjecture}

\numberwithin{equation}{section}

\begin{document}

\title{Quasiperiodic Flows and \\ Algebraic Number Fields}

\maketitle

\markboth{\sc LENNARD F. BAKKER}{\sc QUASIPERIODIC FLOWS}

{\footnotesize
\centerline{LENNARD F. BAKKER}
}
\medskip\medskip\medskip

{\footnotesize
\centerline{ Department of Mathematics, Brigham Young University, Provo, UT 84602 USA }
}

\medskip\medskip\medskip
{\footnotesize
\noindent{\bf ABSTRACT}: We classify a quasiperiodic flow as either algebraic or transcendental. For an algebraic quasiperiodic flow $\phi$ on the $n$-torus, $T^n$, we prove that an absolute invariant of the smooth conjugacy class of $\phi$, known as the multiplier group, is a subgroup of the group of units of the ring of integers in a real algebraic number field $F$ of degree $n$ over ${\mathbb Q}$. We also prove that for any real algebraic number field $F$ of degree $n$ over ${\mathbb Q}$, there exists an algebraic quasiperiodic flow on $T^n$ whose multiplier group is exactly the group of units of the ring of integers in $F$. We support and formulate the conjecture that the multiplier group distinguishes the algebraic quasiperiodic flows from the transcendental ones.
}

\medskip

\section{\bf INTRODUCTION}

\noindent We begin with an example that motivates a connection between a quasiperiodic flow and the group of units of the ring of integers in a real algebraic number field, a connection that manifests itself in terms of an absolute invariant of the smooth conjugacy class of a flow known as the multiplier group.

\begin{example}\label{motivating}
{\rm On the $2$-torus $T^2=S^1\times S^1$ equipped with global coordinates $(\theta_1,\theta_2)$, where $S^1={\mathbb R}/{\mathbb Z}$, let $\phi$ be the quasiperiodic flow generated by the constant vector field
\[ X=\sqrt2\frac{\partial}{\partial\theta_1} 
+\sqrt{10}\frac{\partial}{\partial\theta_2}.\]
The multiplier group of $\phi$ is a subgroup of ${\mathbb R}^*\equiv{\mathbb R}\setminus\{0\}$, the multiplicative group of real numbers, where each element of the multiplier group is an $\alpha\in{\mathbb R}^*$ such that $R_*X=\alpha X$ for a diffeomorphism $R$ of $T^2$. (Here $R_*X={\bf T}RXR^{-1}$ is the push-forward of $X$ by $R$, where ${\bf T}$ is the tangent functor.) For a diffeomorphism $R$ of $T^2$ induced a unimodular matrix $B\in\hbox{\rm GL}(2,{\mathbb Z})$, the equation $R_*X=\alpha X$ becomes the eigenvector-eigenvalue problem $BX=\alpha X$ with a twist: the eigenvector $X$ is known a priori, whereas the matrix $B$ and the eigenvalue $\alpha$ are to be determined. One solution of $BX=\alpha X$ is
\[ \begin{bmatrix} 2 & 1 \\ 5 & 2\end{bmatrix}\begin{bmatrix}\sqrt2 \\ \sqrt{10} \end{bmatrix} = (2+\sqrt 5)\begin{bmatrix}\sqrt2 \\ \sqrt{10}\end{bmatrix}. \]
Two other solutions of $BX=\alpha X$ are $IX=1X$ and $-IX=-1X$ where $I$ is the $2\times 2$ identity matrix. Hence, the set of numbers $G=\{\pm(2+\sqrt 5)^m:m\in{\mathbb Z}\}$ is a subgroup of the multiplier group of $\phi$. On the other hand, $G$ is a proper subgroup of $\{\pm\big((1/2)(1+\sqrt 5)\big)^m:m\in{\mathbb Z}\}$ which is the group of units of the ring of integers in the real algebraic number field ${\mathbb Q}(\sqrt 5)$, where a basis for this field of degree $2$ over ${\mathbb Q}$ is given by the set of components of the scalar multiple $(1/\sqrt 2)X$.
}\end{example}

Guided by Example \ref{motivating}, we formulate two general, but related conjectures about the multiplier group of a quasiperiodic flow $\phi$ on $T^n$ generated by a constant vector field $X$, and the group of units of the ring of integers in a real algebraic number field of degree $n$ over ${\mathbb Q}$.

\begin{conjecture}\label{consubgroup}
If there is a nonzero real $\vartheta$ such that the set of components of the scalar multiple $\vartheta X$ form a basis for a real algebraic number field $F$ of degree $n$ over ${\mathbb Q}$, then the multiplier group of $\phi$ is a subgroup of the group of units of the ring of integers in $F$.
\end{conjecture}

\begin{conjecture}\label{conequality}
If $F$ is a real algebraic number field of degree $n$ over ${\mathbb Q}$, then there is a quasiperiodic flow $\phi$ on $T^n$ whose multiplier group is exactly the group of units of the ring of integers in $F$.
\end{conjecture}

After presenting some generalized symmetry theory (which includes the multiplier group) and some algebraic number theory, we reformulate and prove these conjectures (see Theorem \ref{subgroup} and Theorem \ref{equality}). As a consequence of the reformulation, we define the algebraic and transcendental classes into which the set of all quasiperiodic flows separates (see Definition \ref{algtrans}), and support and formulate the conjecture that the multiplier group distinguishes between these two classes.

\section{\bf MULTIPLIERS AND QUASIPERIODIC FLOWS}

\noindent For $n\geq2$ an integer, let $T^n=S^1\times\cdot\cdot\cdot\times S^1$ ($n$-times) be the $n$-torus which is a smooth, i.e. $C^\infty$, manifold without boundary, and let $\hbox{\rm Diff}(T^n)$ be the group of smooth diffeomorphisms of $T^n$. Let $\phi:{\mathbb R}\times T^n\to T^n$ be a smooth (complete) flow on $T^n$, and let
\[ X(p)=\frac{d}{dt}\phi(t,p)\bigg\vert_{t=0}\]
be the smooth vector field that generates $\phi$. A flow $\phi$ on $T^n$ is {\it nontrivial} if its generating vector field $X$ is not identically zero.

The generalized symmetries of a flow considered here are a type of space-time symmetry that sends orbits to orbits in which a global, uniform, nonzero constant time reparameterization of all orbits is permitted. This type of space-time symmetry can be simply defined in terms of the push-forward of a diffeomorphism on the generating vector field $X$ of a flow and a nonzero scalar multiple of $X$.

\begin{definition}\label{gensym}
An $R\in\hbox{\rm Diff}(T^n)$ is a generalized symmetry of $\phi$ if there is an $\alpha\in{\mathbb R}^*$ such that $R_*X=\alpha X$. The number $\alpha$ is a multiplier of $R$ for $\phi$.
\end{definition}

Let $S_\phi$ be the set of the generalized symmetries of $\phi$ in the sense of Definition \ref{gensym}. This set is a subgroup of $\hbox{\rm Diff}(T^n)$ (see Proposition 2.1 in \cite{CGSF}); it is called the {\it generalized symmetry group} of $\phi$; and it is the group theoretic normalizer of $F_\phi=\{\phi_t:t\in{\mathbb R}\}$, the abelian group of diffeomorphisms determined by the flow where $\phi_t$ is defined by $\phi_t(\theta)=\phi(t,\theta)$ (Theorem 2.5 in \cite{CGSF}). The generalized symmetry group of a flow $\phi$ is a proper subgroup of $\hbox{\rm Diff}(T^n)$ if and only if $\phi$ is nontrivial (see Theorem 2.6 in \cite{CGSF}). In Definition \ref{gensym}, the value of a multiplier of a generalized symmetry was not assumed to be unique.

\begin{lemma}
For a nontrivial flow $\phi$, each $R\in S_\phi$ has a unique multiplier.
\end{lemma}

\begin{proof}
See Lemma 2.7 in \cite{CGSF}.
\end{proof}

For a nontrivial flow $\phi$, the {\it multiplier representation}\, $\rho_\phi:S_\phi\to{\mathbb R}^*\cong\hbox{\rm GL}({\mathbb R})$ is a linear representation of $S_\phi$ in ${\mathbb R}$ (Theorem 2.8 in \cite{CGSF}), a homomorphism which takes each $R\in S_\phi$ to its unique multiplier $\rho_\phi(R)$. The image of this multiplier representation, $\rho_\phi(S_\phi)$, is the {\it multiplier group}\, of $\phi$. The multiplier group of a nontrivial flow $\phi$ always contains $1$: $R=\hbox{\rm id}$, the identity diffeomorphism of $T^n$, and $\alpha=1$ satisfy $R_*X=\alpha X$.

Two flows $\phi$ and $\psi$ on $T^n$ generated by $X$ and $Y$ respectively are {\it smoothly conjugate} if and only if there is a $V\in\hbox{\rm Diff}(T^n)$ such that $V_*X=Y$. Smooth conjugacy is an equivalence relation on the set of smooth flows on $T^n$, and hence divides the set of smooth flows on $T^n$ into smooth conjugacy classes. An {\it absolute invariant} of the smooth conjugacy class of a nontrivial $\phi$ is its multiplier group.

\begin{theorem}
If two nontrivial flows $\phi$ and $\psi$ on $T^n$ are smoothly conjugate, then the multiplier groups $\rho_\phi(S_\phi)$ and $\rho_\psi(S_\psi)$ are identical $($and not just isomorphic$)$.
\end{theorem}

\begin{proof}
See Theorem 4.2 in \cite{CGSF}.
\end{proof}

A flow $\phi$ on $T^n$ generated by $X$ is {\it quasiperiodic} if and only if it is smoothly conjugate with a flow on $T^n$ generated by a constant vector field whose $n$ real components are independent over ${\mathbb Q}$ (pp.\ 79-80 \cite{HB}). (Here $n$ real numbers $a_1,...,a_n$ are {\it independent over} ${\mathbb Q}$, or rationally independent, if for $(m_1,...,m_n)\in{\mathbb Z}^n$, the equation $\sum_{j=1}^n m_j a_j=0$ implies that $m_j=0$ for all $j=1,...,n$.) A quasiperiodic flow is necessarily nontrivial, and its multiplier group always contains $-1$ in addition to $1$ (Theorem 2.1 in \cite{SGIQ}). Since the multiplier group is an absolute invariant of a nontrivial flow's smooth conjugacy class, and since, by definition, any quasiperiodic flow is smoothly conjugate to a quasiperiodic flow generated by a constant vector field, attention is restricted to quasiperiodic flows on $T^n$ that are generated by constant vector fields.

In terms of fixed, global coordinates $(\theta_1,...,\theta_n)$ on $T^n$, a constant vector field $X$ on $T^n$ has the form $X=\sum_{i=1}^n a_i\partial/\partial\theta_i$ where $a_i\in{\mathbb R}$. The equation $R_*X=\alpha X$ that defines a generalized symmetry becomes ${\bf T}RX=\alpha X$, and in this form it restricts the possibilities for the multipliers of a quasiperiodic flow.

\begin{theorem}\label{unimodular}
Let $\phi$ be a quasiperiodic flow on $T^n$ generated by a constant vector field $X=\sum_{i=1}^n a_i\partial/\partial\theta_i$. If $\alpha\in\rho_\phi(S_\phi)$, then there exists a unique $B=(b_{ij})\in\hbox{\rm GL}(n,{\mathbb Z})$ such that
\[ \alpha a_i = \sum_{j=1}^n b_{ij}a_j{\rm\ for\ each\ }i=1,...,n.\]
\end{theorem}

\begin{proof}
See Corollary 4.5 and Theorem 5.3 in \cite{SGIQ}.
\end{proof}

The $n$ equations in Theorem \ref{unimodular} that the entries of the unimodular matrix $B$ and the multiplier $\alpha$ satisfy are the components of ${\bf T}RX=\alpha X$ relative to the global coordinates $(\theta_1,...,\theta_n)$. The $n$ equations are derived from the result that the tangent map ${\bf T}R$\, for every $R\in S_\phi$ belongs to $\hbox{\rm GL}(n,{\mathbb Z})$ (Theorem 4.3 in \cite{SGIQ}). Thus, the search for the generalized symmetries and the concomitant multipliers of a quasiperiodic flow reduces to the eigenvector-eigenvalue problem $BX=\alpha X$ with a twist as shown in Example \ref{motivating}: the generating vector field $X$ is the given eigenvector, while a solution is a pair $(B,\alpha)$ that satisfies $BX=\alpha X$, in which the unimodular matrix $B$ is uniquely determined by the eigenvalue $\alpha$ as stated in Theorem \ref{unimodular}. A complete characterization of the relationship between the components of $X$ and the entries of $B$ is known for any quasiperiodic flow on $T^2$ (Theorem 2.3 in \cite{RRSGQF}) and was used in Example \ref{motivating} to find the $B$ given there; such a characterization is also known for a certain kind of quasiperiodic flow on $T^3$ (Theorem 3.1 in \cite{RRSGQF}).

\begin{example}\label{2torus1}{\rm 
Let $\phi$ be the quasiperiodic flow on $T^2$ generated by $X=\partial/\partial\theta_1+(1+\sqrt2)\partial/\partial\theta_2$. Then $R\in S_\phi$ if and only if there exist unique $u_1,u_2\in{\mathbb Z}$ such that
\[ B={\bf T}R=\begin{bmatrix}-2u_1+u_2 & u_1 \\ u_1 & u_2 \end{bmatrix} \]
with $\hbox{\rm det}(B)=\pm 1$ (Theorem 2.3 in \cite{RRSGQF}). The value of $\rho_\phi(R)$ is determined by any of the $n$ equations listed in the conclusion of Theorem \ref{unimodular}; the $i=1$ equation gives
\[\rho_\phi(R)=(-2u_1+u_2)+(1+\sqrt 2)u_1. \]
The choice of $u_1=1$ and $u_2=2$ gives
\[ B=\begin{bmatrix} 0 & 1 \\ 1 & 2\end{bmatrix} \]
whose determinant is $-1$, so that $B$ is unimodular and the diffeomorphism induced by $B$ is a generalized symmetry of $\phi$ whose multiplier is $1+\sqrt2$.
}\end{example}

\begin{example}\label{3torus1}{\rm
Let $\beta$ be the real root of the irreducible $z^3+z-1$. (An exact expression for $\beta$ is given by Cardano's cubic formula.) Let $\psi$ be the quasiperiodic flow on $T^3$ generated by $Y=\beta^2\partial/\partial\theta_1+\beta\partial/\partial\theta_2+\partial/\partial\theta_3$. Then $R\in S_\psi$ if and only if there are unique $u_1,u_2,u_3\in{\mathbb Z}$ such that
\[ B={\bf T}R=\begin{bmatrix}u_1+u_2+u_3 & u_1 & u_1+u_2 \\ u_1+u_2 & u_3 & u_1 \\ u_1 & u_2 & u_3\end{bmatrix}\]
with $\hbox{\rm det}(B)=\pm1$ (Theorem 3.1 in \cite{RRSGQF}). The value of $\rho_\psi(R)$ is determined by the $i=3$ equation listed in the conclusion of Theorem \ref{unimodular}:
\[ \rho_\psi(R)=u_1\beta^2+u_2\beta+u_3.\]
The choice of $u_1=u_2=u_3=1$ gives
\[ B=\begin{bmatrix} 3 & 1 & 2 \\ 2 & 1 & 1 \\ 1 & 1 & 1 \end{bmatrix} \]
whose determinant is $1$, so that $B\in\hbox{\rm GL}(3,{\mathbb Z})$ and the diffeomorphism of $T^3$ induced by $B$ is a generalized symmetry of $\psi$ whose multiplier is $\beta^2+\beta+1$.
}\end{example}

\section{\bf MULTIPLIERS AND ALGEBRAIC NUMBER FIELDS}

\noindent An {\it algebraic number field} $F$ is a finite dimensional field extension of ${\mathbb Q}$. The {\it degree of} $F$ over ${\mathbb Q}$ is the dimension of $F$ as a vector space over ${\mathbb Q}$. A {\it basis} for a real algebraic number field $F$ of degree $n$ over ${\mathbb Q}$ is a subset $\{a_1,...,a_n\}$ of $F$ whose elements are independent over ${\mathbb Q}$ and span $F$; that is, every $\alpha\in F$ is uniquely expressed as $\alpha=\sum_{i=1}^n \beta_i a_i$ where $\beta_i\in{\mathbb Q}$.

\begin{definition}\label{Falgebraic}
Let $F$ be a real algebraic number field of degree $n$ over ${\mathbb Q}$. A quasiperiodic flow $\phi$ on $T^n$ generated by the constant vector field $X=\sum_{i=1}^n a_i\partial/\partial\theta_i$ is called $F$-algebraic if there is a $\vartheta\in{\mathbb R}^*$ such that $\{\vartheta a_1,...,\vartheta a_n\}$ is a basis for $F$.
\end{definition}

As an invariant of smooth conjugacy, the multiplier group is not complete: two flows in different smooth conjugacy classes can have the same multiplier group (see Example 4.6 in \cite{CGSF}). This incompleteness manifests itself for two quasiperiodic flows that are both $F$-algebraic when the generating vector field of one flow is a nonzero scalar multiple of the generating vector field of the other flow.

\begin{proposition}\label{Falgident}
Let $F$ be an algebraic number field of degree $n$ over ${\mathbb Q}$, and let $\phi$ and $\psi$ be quasiperiodic flows generated by the constant vector fields $X=\sum_{i=1}^n a_i\partial/\partial\theta_i$ and $Y=\sum_{i=1}^n\omega_i\partial/\partial_i$ respectively. If $\phi$ is $F$-algebraic and $X=\mu Y$ for $\mu\in{\mathbb R}^*$, then $\rho_\phi(S_\phi)=\rho_\psi(S_\psi)$.
\end{proposition}

\begin{proof}
Under the hypotheses, there are $\vartheta,\mu\in{\mathbb R}^*$ such that the set of components of $\vartheta X$ is a basis for $F$, and $X=\mu Y$. Let $\alpha\in\rho_\phi(S_\phi)$. Then by Theorem \ref{unimodular} and the discussion following its proof, there is a unique $B\in\hbox{\rm GL}(n,{\mathbb Z})$ such $BX=\alpha X$. The equation $BX=\alpha X$ is $R_*X=\alpha X$ where $R$ is the diffeomorphism of $T^n$ induced by $B$. Since $X=\mu Y$, the equation $R_*X=\alpha X$ implies the equation $R_*Y=\alpha Y$. So $\alpha\in\rho_\psi(S_\psi)$. Reversing the roles of $X$ and $Y$ gives the other inclusion.
\end{proof}

A real number is an {\it algebraic integer} if it is a real root of a monic polynomial with coefficients in ${\mathbb Z}$ (see pp.\ 1-2 in \cite{S-D}). The {\it degree} of an algebraic integer is the degree of the irreducible monic polynomial with coefficients in ${\mathbb Z}$ for which the algebraic integer is a root (see pg.\ 1 in \cite{MO}).

\begin{lemma}\label{integer}
Let $\phi$ be a quasiperiodic flow on $T^n$ generated by a constant vector field. 
If $\alpha\in\rho_\phi(S_\phi)$, then $\alpha$ is a real algebraic integer of degree at most $n$.
\end{lemma}

\begin{proof}
See Corollary 4.4 in \cite{SGIQ}.
\end{proof}

Let $\alpha$ be an element of a real algebraic number field $F$ of degree $n$ over ${\mathbb Q}$, and choose a basis $\{a_1,...,a_n\}$ for $F$. The {\it norm} of $\alpha$, denoted by $\hbox{\rm norm}_{F/{\mathbb Q}}(\alpha)$, is the determinant of the unique $n\times n$ matrix $B=(b_{ij})$ with $b_{ij}\in{\mathbb Q}$ such that
\[ \alpha a_i = \sum_{j=1}^n b_{ij}a_j\quad i=1,...,n\]
(see p.\ 121 in \cite{S-D}). The value of $\hbox{\rm norm}_{F/{\mathbb Q}}(\alpha)$ is independent of the basis chosen for $F$. The similarity of the definition of the norm of an element of $F$ and the conclusion of Theorem \ref{unimodular} is the key to proving Conjecture \ref{consubgroup}.

Let ${\mathfrak o}_F$ denote the {\it ring of integers} in a real algebraic number field $F$, that is the set of real algebraic integers in $F$. The degree of an  $\alpha\in{\mathfrak o}_ F$ does not exceed the degree $n$ of $F$ over ${\mathbb Q}$ because every element of $F$ is a real root of a monic polynomial with coefficients in ${\mathbb Q}$ of degree at most $n$ (see Theorem 4.9, p.\ 50 and Theorem 5.6, p.\ 62 in \cite{P-D}). For a real algebraic number field $F$ of degree $n$ over ${\mathbb Q}$, there is an $\alpha\in{\mathfrak o}_F$ of degree $n$ such that $F$ is the smallest subfield of ${\mathbb R}$ containing ${\mathbb Q}\cup\{\alpha\}$ (see Theorem 5.6, p.\ 62 and Theorem 6.5 p.\ 77 in \cite{P-D}); this smallest subfield is denoted by ${\mathbb Q}(\alpha)$.

A {\it unit} of ${\mathfrak o}_F$ is an $\alpha\in{\mathfrak o}_F$ whose multiplicative inverse, $\alpha^{-1}$, is also in ${\mathfrak o}_F$. The {\it group of units}, denoted by ${\mathfrak o}^*_F$, is the set of units in ${\mathfrak o}_F$. A useful characterization of a unit is: $\alpha\in{\mathfrak o}^*_F$ if and only if $\alpha\in{\mathfrak o}_F$ and $\hbox{\rm norm}_{F/{\mathbb Q}}(\alpha)=\pm1$ (see p.\ 20 in \cite{S-D}).

\begin{theorem}\label{subgroup}
Let $\phi$ be a quasiperiodic flow on $T^n$ generated by the constant vector field $X=\sum_{i=1}^n a_i\partial/\partial\theta_i$. If $\phi$ is $F$-algebraic, then $\rho_\phi(S_\phi)<{\mathfrak o}_F^*$. 
\end{theorem}

\begin{proof}
By the hypotheses, there is a nonzero real $\vartheta$ such that $\{\vartheta a_1,...,\vartheta a_n\}$ is a basis for the real algebraic number field $F$ of degree $n$ over ${\mathbb Q}$. Let $\alpha\in\rho_\phi(S_\phi)$. By Theorem \ref{unimodular} there exists a unique $B=(b_{ij})\in\hbox{\rm GL}(n,{\mathbb Z})$ such that
\[ \alpha=\sum_{j=1}^n b_{ij}\frac{a_j}{a_i} = \frac {1}{\vartheta a_i}\sum_{j=1}^n b_{ij}\vartheta a_j\]
for every $i=1,...,n$. It follows that $\alpha\in F$ since $F$ is a field and $\{\vartheta a_1,...,\vartheta a_n\}$ is a basis for $F$. Moreover, $\alpha\in{\mathfrak o}_F$ by Lemma \ref{integer}. Finally, $\alpha\in{\mathfrak o}_F^*$ since $\hbox{\rm norm}_{F/{\mathbb Q}}(\alpha)=\hbox{\rm det}(B)=\pm 1$.
\end{proof}

\begin{example}\label{2torus2}{\rm 
Let $\psi$ be the quasiperiodic flow on $T^2$ generated by $Y=\partial/\partial\theta_1+\sqrt 5\partial/\partial\theta_2$. This flow is $F$-algebraic for $F={\mathbb Q}(\sqrt5)$, and the generating vector field $X=\sqrt2\partial/\partial\theta_1+\sqrt{10}\partial/\partial\theta_2$ of the flow $\phi$ in Example \ref{motivating} is a nonzero scalar multiple of $Y$. By Proposition \ref{Falgident} the multiplier groups of $\psi$ and $\phi$ are identical. By Theorem \ref{subgroup} this common multiplier group is a subgroup of ${\mathfrak o}^*_F$. However, $\rho_\psi(S_\psi)\ne{\mathfrak o}^*_F$ because the unit $(1/2)(1+\sqrt 5)$ (a root of $z^2-z-1$) is excluded from $\rho_\psi(S_\psi)$ by Theorem \ref{unimodular}: it is impossible that 
$(1/2)(1+\sqrt5)=b_{11}+b_{12}\sqrt5$ for $b_{11},b_{12}\in{\mathbb Z}$. Also, $\phi$ and $\psi$ belong to different smooth conjugacy classes since $V_*X=Y$ requires a diffeomorphism $V$ of $T^2$ such that ${\bf T}V=\sqrt2 I$, which is impossible.
}\end{example}

\begin{example}\label{4torus1}{\rm
The real algebraic number field $F={\mathbb Q}(\sqrt2+\sqrt3)$ has degree $4$ over ${\mathbb Q}$ because $\sqrt2+\sqrt3$ is a root of the irreducible $z^4-10z^2+1$. A basis for $F$ is $\{1,\sqrt2,\sqrt3,\sqrt6\}$. Let $\phi$ be the quasiperiodic flow on $T^4$ generated by
\[ X=\sqrt6\frac{\partial}{\partial\theta_1}
+\sqrt3\frac{\partial}{\partial\theta_2}
+\sqrt2\frac{\partial}{\partial\theta_2}
+\frac{\partial}{\partial\theta_4}.\]
This quasiperiodic flow is thus $F$-algebraic, so that $\rho_\phi(S_\phi)$ is a subgroup of ${\mathfrak o}_F^*$ by Theorem \ref{subgroup}. In particular, the number $\alpha=2\sqrt6+4\sqrt3+5\sqrt2+5$ is a multiplier of $\phi$ because the diffeomorphism $R$ of $T^4$ induced by the unimodular matrix
\[ B=\begin{bmatrix} 5 & 10 & 12 & 12 \\ 5 & 5 & 6 & 12 \\ 4 & 4 & 5 & 10 \\ 2 & 4 & 5 & 5\end{bmatrix} \]
satisfies $R_*X=\alpha X$. Hence $\alpha\in{\mathfrak o}_F$. However, $\rho_\phi(S_\phi)\ne{\mathfrak o}_F^*$ because the unit $(1/2)(\sqrt2+\sqrt6)$ (a root of $z^4-4z^2+1$) is excluded from $\rho_\phi(S_\phi)$ by Theorem \ref{unimodular}: it is impossible that
\[\frac{\sqrt2+\sqrt6}{2}=c_{41}\sqrt6+c_{42}\sqrt3+c_{43}\sqrt2+c_{44}\]
for integers $c_{4j}$, $j=1,2,3,4$.
}\end{example}

The ring of integers ${\mathfrak o}_F$ is a ${\mathbb Z}$-module (see Lemma 1, p. 2 in \cite{S-D}). An {\it integral basis} for a real algebraic number field $F$ of degree $n$ over ${\mathbb Q}$ is a set of real algebraic integers $\{a_1,...,a_n\}$ in $F$ that is a ${\mathbb Z}$-module basis for ${\mathfrak o}_F$ and is also a basis for $F$ (see pp.\ 79-80 in \cite{P-D}). The existence of an integral basis for $F$ is the key to proving Conjecture \ref{conequality}.

\begin{lemma}\label{intbasisexist}
If $F$ is a real algebraic number field, then there exists an integral basis for $F$.
\end{lemma}

\begin{proof}
See Theorem 6.9, p. 81 in \cite{P-D}.
\end{proof}

\begin{theorem}\label{equality}
If $F$ is a real algebraic number field of degree $n$ over\, ${\mathbb Q}$, then there exists an $F$-algebraic quasiperiodic flow $\phi$ on $T^n$ such that $\rho_\phi(S_\phi)={\mathfrak o}_F^*$.
\end{theorem}

\begin{proof}
Suppose $F$ is a real algebraic number field of degree $n$ over ${\mathbb Q}$. By Lemma \ref{intbasisexist} there is an integral basis for $F$, say $\{a_1,...,a_n\}$. Let $\phi$ be the flow on $T^n$ generated by $X=\sum_{i=1}^n a_i\partial/\partial\theta_i$. Since the $a_1,...,a_n$ are independent over ${\mathbb Q}$, the flow $\phi$ is quasiperiodic. Since $\{a_1,...,a_n\}$ is a basis for $F$, the flow $\phi$ is $F$-algebraic. Then $\rho_\phi(S_\phi)\subset{\mathfrak o}_F^*$ by Theorem \ref{subgroup}, so it suffices to show that ${\mathfrak o}_F^*\subset\rho_\phi(S_\phi)$. Let $\alpha\in{\mathfrak o}_F^*$. Then since $\hbox{\rm norm}_{F/{\mathbb Q}}(\alpha)=1$ and since $\{a_1,...,a_n\}$ is an integral basis for $F$, there exists a unique $B=(b_{ij})\in\hbox{\rm GL}(n,{\mathbb Z})$ such that
\[ \alpha a_i = \sum_{j=1}^n b_{ij}a_j\]
for $i=1,...,n$. Let $R$ be the diffeomorphism of $T^n$ induced by the unimodular matrix $B$. Then the equations $\alpha a_i=\sum_{j=1}^n b_{ij} a_j$, $i=1,...,n$, are the components of $\alpha X=R_*X$. Hence $R\in S_\phi$, and $\alpha=\rho_\phi(R)\in\rho_\phi(S_\phi)$.
\end{proof}

A quasiperiodic flow on $T^n$ whose multiplier group is the group of units of the ring of integers in a real algebraic number field $F$ of degree $n$ over ${\mathbb Q}$ is constructed by finding an integral basis for $F$ and using the elements of this integral basis as the components of the generating vector field.

\begin{example}\label{2torus3}{\rm 
An integral basis for $F={\mathbb Q}(\sqrt2)$ is $\{1,\sqrt2\}$ (see Theorem 6.11, p.\ 83 in \cite{P-D}). The multiplier group of the quasiperiodic flow $\psi$ on $T^2$ generated by $Y=\partial/\partial\theta_1+\sqrt2\partial/\partial\theta_2$ is ${\mathfrak o}^*_F$. This group is $\{\pm(1+\sqrt2)^m:m\in{\mathbb Z}\}$ (Theorem 7.9, p.\ 92 in \cite{P-D}), which is isomorphic to ${\mathbb Z}_2\times{\mathbb Z}$. The multiplier group of the quasiperiodic flow $\phi$ in Example \ref{2torus1}, generated by $X=\partial/\partial\theta_1+(1+\sqrt2)\partial/\partial\theta_2$, is also ${\mathfrak o}^*_F$ because $\{1,1+\sqrt2\}$ is an integral basis for $F$ (the number $1+\sqrt2$ is a root of $z^2-2z-1$). The flows $\phi$ and $\psi$ lie in the same smooth conjugacy class because $V_*X=Y$ for the diffeomorphism $V$ of $T^2$ induced by the unimodular matrix
\[\begin{bmatrix}1 & 0 \\ -1 & 1\end{bmatrix}.\]
This is typical because different integral bases of the same algebraic number field are related by a unimodular matrix.
}\end{example}

\begin{example}\label{2torus4}{\rm
An integral basis for $F={\mathbb Q}(\sqrt 5)$ is $\{1,(1/2)(1+\sqrt 5)\}$ (see Theorem 6.11, p.\ 83 in \cite{P-D}). The multiplier group of the quasiperiodic flow $\chi$ on $T^2$ generated by $Z=\partial/\partial\theta_1+(1/2)(1+\sqrt 5)\partial/\partial\theta_2$ is ${\mathfrak o}^*_F$. This group is $\{\pm\big((1/2)(1+\sqrt 5)\big)^m:m\in{\mathbb Z}\}$, which is isomorphic to ${\mathbb Z}_2\times{\mathbb Z}$.
}\end{example}

\begin{example}\label{4torus2}{\rm
An integral basis for $F={\mathbb Q}(\sqrt2+\sqrt3)$ is $\{1,\sqrt3,\sqrt 6,(1/2)(\sqrt 2+\sqrt6)\}$ (see pp.\ 84-85 in \cite{HC}). The multiplier group of the flow $\psi$ generated by
\[ Y=\frac{\sqrt2+\sqrt6}{2}\frac{\partial}{\partial\theta_1}
+\sqrt6\frac{\partial}{\partial\theta_2}
+\sqrt3\frac{\partial}{\partial\theta_3}
+\frac{\partial}{\partial\theta_4}\]
is ${\mathfrak o}^*_F$. This group is isomorphic to ${\mathbb Z}_2\times{\mathbb Z}\times{\mathbb Z}\times{\mathbb Z}$ by Dirichlet's Unit Theorem (p.\ 21 in \cite{S-D}) because $\sqrt2+\sqrt3$ is a root of the irreducible $z^4-10z^2+1$ which has four real roots.
}\end{example}

\section{\bf THE ALGEBRAIC AND TRANSCENDENTAL CLASSES}

\noindent For any real algebraic number field $F$, the abelian group ${\mathfrak o}^*_F$ is finitely generated and contains at least one infinite cyclic factor by Dirichlet's Unit Theorem. The multiplier groups of the quasiperiodic flows in the examples given thus far each contain at least one infinite factor because there exists a multiplier other than $\pm1$. This coincidence leads to the formulation of two classes of quasiperiodic flows and a conjecture regarding their separation in terms of the multiplier group.

\begin{definition}\label{algtrans}
A quasiperiodic flow $\phi$ on $T^n$ generated by a constant vector field $X$ is algebraic if it is $F$-algebraic for a real algebraic number field $F$ of degree $n$ over ${\mathbb Q}$; otherwise $\phi$ is transcendental.
\end{definition}

The multiplier group for the algebraic quasiperiodic flow $\phi$ in Example \ref{motivating} is the group $G=\{\pm(2+\sqrt5)^m:m\in{\mathbb Z}\}$ because, $\rho_\phi(S_\phi)$ being a subgroup of ${\mathfrak o}^*_F=\{\pm\big((1/2)(1+\sqrt 5)\big)^m:m\in{\mathbb Z}\}$ for $F={\mathbb Q}(\sqrt 5)$ by Theorem \ref{subgroup}, the multiplier $2+\sqrt5$ of $\phi$ is the smallest unit in ${\mathfrak o}^*_F$ that is greater than $1$ and in $\rho_\phi(S_\phi)$. The index of $\rho_\phi(S_\phi)$ as a subgroup of ${\mathfrak o}^*_F$ is $3$ because
$(1/2)(1+\sqrt 5)\not\in G$, $\big((1/2)(1+\sqrt 5)\big)^2\not\in G$, while
$\big((1/2)(1+\sqrt5)\big)^3=2+\sqrt5\in G$.

The multiplier group of the quasiperiodic flow $\psi$ on $T^2$ generated by $Y=\partial/\partial\theta_1+\pi\partial/\partial\theta_2$ is $\{1,-1\}$ (Theorem 2.1 in \cite{RRSGQF}). Suppose that $\psi$ is algebraic. Then there is a real algebraic number field $F$ of degree $2$ over ${\mathbb Q}$ and a $\vartheta\in{\mathbb R}^*$ such that $\{\vartheta,\vartheta\pi\}$ is a basis for $F$. Thus $\pi=\vartheta\pi/\vartheta\in F$ since $F$ is a field. This implies that $\pi$ is a root of a monic polynomial with coefficients in ${\mathbb Q}$ of degree at most $2$. This contradiction shows that $\psi$ is transcendental.

\begin{conjecture}
The multiplier group distinguishes the algebraic quasiperiodic flows from the transcendental ones in that $($a$)$ $\phi$ is an algebraic quasiperiodic flow on $T^n$ if and only if $\rho_\phi(S_\phi)$ is a finite index subgroup of ${\mathfrak o}^*_F$ for a real algebraic number field $F$ of degree $n$ over ${\mathbb Q}$, and $($b$)$ $\phi$ is a transcendental quasiperiodic flow on $T^n$ if and only if $\rho_\phi(S_\phi)=\{1,-1\}$.
\end{conjecture}

\end{document}